\documentclass{article}
\pdfoutput=1
\usepackage{amsmath}
\usepackage{cases}
\usepackage{amsfonts}
\usepackage{array}
\usepackage{geometry}
\title{\textbf{On a New Formula for Arithmetic Functions}}
\author{J.AKOUN}

\begin{document}
	\maketitle	
	
	\begin{abstract}
\centering\begin{minipage}{\dimexpr\paperwidth-10cm}
In this paper we establish a new formula for the arithmetic functions that verify $ f(n) = \sum_{d|n} g(d)$ where $g$ is also an arithmetic function. We prove the following identity,
$$\forall n \in \mathbb{N}^*, \ \ \ f(n) = \sum_{k=1}^n \mu \left(\frac{k}{(n,k)}\right)  \frac {\varphi(k)}{\varphi\left(\frac{k}{(n,k)}\right)}   \sum_{l=1}^{\left\lfloor\frac{n}{k}\right\rfloor} \frac{g(kl)}{kl} $$
where $\varphi$ and  $\mu$ are respectively Euler's and Mobius' functions and (.,.) is the  GCD. First, we will compare this expression with other known expressions for arithmetic functions and pinpoint its advantages. Then, we will prove the identity using exponential sums' proprieties. Finally we will present some applications with well known functions such as $d$ and $\sigma$ which are respectively the number of divisors function and the sum of divisors function.
\end{minipage}
\end{abstract}
\bigskip

	\section{Introduction}
\subparagraph{}	
In the paper we present and prove a new transformed expression for some arithmetic functions. First let's remind the definition of an arithmetic function,\\
\\
\textbf{Definition.}
\textit{A function $f$ is arithmetic if its domain is the positive integers and hence \ \  $ f : \mathbb{N}^* \longrightarrow \mathbb{C}$.}\\
\\
The main theorem is,\\
\\
\textbf{Theorem.}
\textit{If $ f(n) = \sum_{d|n} g(d)$ where $g$ is an arithmetic function then,
$$\forall n \in \mathbb{N}^*, \ \ \ f(n) = \sum_{k=1}^n \mu \left(\frac{k}{(n,k)}\right)  \frac {\varphi(k)}{\varphi\left(\frac{k}{(n,k)}\right)}   \sum_{l=1}^{\left\lfloor\frac{n}{k}\right\rfloor} \frac{g(kl)}{kl} $$
where $\varphi$ and  $\mu$  are respectively Euler's and Mobius' functions and (.,.) is the GCD.}\\
\\
The strength of this new expression is that the sum is not indexed on divisors but on all the integers. As a matter of fact let's take the example of $\sigma$ which is the sum of divisors function. We can clearly write $\sigma$ as a sum with,
$$\sigma(n) = \sum_{k|n} k $$
however this formula is not really interesting because we don't control the indexation. With this expression we  moved the unknown from $\sigma$ to the indexation, whereas our transformed expression for arithmetic functions does not hide complexity in the indexation. For $\sigma$ the theorem gives us,

$$\forall n \in \mathbb{N}^*, \ \ \  \sigma(n) = \sum_{k=1}^n \left\lfloor\frac{n}{k}\right\rfloor  \mu\left(\frac{k}{(n,k)}\right)  \frac {\varphi(k)}{\varphi\left(\frac{k}{(n,k)}\right)} $$

Still, the uncontrolled arithmetical part has not completely vanished and lies in Euler's and Mobius' functions. 
\subparagraph{}
The idea to express arithmetical functions in terms of sums was first explored by Ramanujan's who gave expressions with series. For example he proved that,

$$\forall n \in \mathbb{N}^*, \ \ \  \sigma(n) =\frac{\pi^2}{6} \sum_{k=1}^{\infty} \frac{n}{k^2} \mu\left(\frac{k}{(n,k)}\right)  \frac {\varphi(k)}{\varphi\left(\frac{k}{(n,k)}\right)} $$
\\
Our expression is reminiscent of Ramanujan's one. This will turn out coherent since the proof of our theorem relies on Ramanujan's sums. Nevertheless let's underline that our theorem gives expressions with finite sums that might be easier to manipulate than series.
 
\subparagraph{} 
Finally, finding a new formulas for arithmetical functions is always exciting as they are at the core of great modern problems. For instance, $\sigma$ is closely linked to the Riemann hypothesis (RH) as underline the two following equivalences.\\
\\
\textit{\textbf{Robin's equivalence}}
$$RH\ \  \Longleftrightarrow \ \ \forall n \geq 5041, \ \ \ \sigma(n) < ne^{\gamma}log(log(n))$$
where $\gamma$ is Euler-Mascheroni constant.\\
\\
\textit{\textbf{Lagarias' equivalence}}
$$RH\ \  \Longleftrightarrow \ \ \forall n > 1, \ \ \ \sigma(n) < H_n+log(H_n)e^{H_n}$$
\\
where $H_n = \sum_{k=1}^n \frac{1}{k}$.

	\section{Preliminary results on exponential sums}	
\subsection{Classical exponential sums} 
\subparagraph{}
First we will focus on classical exponential sums. Those sums are well known and will be useful to establish the theorem presented in this article. For those sums we know that,\\
\begin{equation}
\forall n \in \mathbb{N}^*, \ \forall m \in \mathbb{N},\ \ \ \sum_{k=1}^n \exp\left(\frac{2i\pi km}{n}\right) =
    \begin{cases}
      n & if\  n | m \\
      0 & otherwise
    \end{cases}
\end{equation}
\\
This result can easily be established with the formula of geometric sums.
\subsection{Ramanujan's sums}
\subparagraph{}
The second type of exponential sums which will turn out useful for the proof are Ramanujan's sums. They are defined by,
$$ \mathbf{c}_m(n) = \sum_{1 \leq k \leq n \atop (k,m) = 1} \exp\left(\frac{2i\pi kn}{m}\right) $$
Those sums are harder to study than the ones before. Still we can remark that,
\subparagraph{•}
If $(m,n)=1$, 

\begin{equation*}
\Phi =
    \begin{cases}
    \mathbb{Z} / m \mathbb{Z} \rightarrow \mathbb{Z} / m \mathbb{Z} & \\
      x \longmapsto n \cdot x &
    \end{cases}
\end{equation*} 
\\
is an isomorphism. Thus,
 
$$\mathbf{c}_m(n) = \sum_{1 \leq k \leq m \atop (k,m) = 1} \exp\left(\frac{2i\pi kn}{m}\right) = \sum_{1 \leq k \leq m \atop (k,m) = 1} \exp\left(\frac{2i\pi k}{m}\right) = \mathbf{c}_m(1)$$
\\
moreover we know that,
\begin{equation}
\sum_{d | n} \mu(d) =
    \begin{cases}
    1 & if\ \ n = 1 \\
    0 & otherwise
    \end{cases}
\end{equation}
if we rearrange the sum it comes,
\\
\begin{align*}
\forall m \in \mathbb{N}^*, \ \ \ 
\mathbf{c}_m(1) &= \sum_{1 \leq k \leq m \atop (k,m) = 1} \exp\left(\frac{2i\pi k}{m}\right)  \\
&= \ \ \sum_{k=1}^m \exp\left(\frac{2i\pi k}{m}\right) \sum_{d | (k,m)} \mu(d) \\
&=\ \ \sum_{d|m} \mu(d) \sum_{l=1}^\frac{m}{d} \exp\left(\frac{2i\pi l}{m}\right)\\
\intertext{with the equality (1) it comes,}\\
&=\ \ \sum_{d|m} \mu(d) \cdot \mathbf{1}_{d=m} \\
&= \ \ \mu(m)
\end{align*}
 
\subparagraph{•}
If $(m,n) =n$, 
$$\mathbf{c}_m(n) = \sum_{1 \leq k \leq n \atop (k,m) = 1} \exp\left(\frac{2i\pi km}{n}\right) = \sum_{1 \leq k \leq n \atop (k,m) = 1} 1 = \varphi(m)$$
\\
\\
\textbf{Remark.} We already see with those two examples that Ramanujan's sums are closely linked with Euler's and Mobius' functions.\\
\\
Except for those values of $(m,n)$ the expression of $\mathbf{c}_m(n)$ is harder to find. Fortunately Hölder showed in 1936 that,\\
\begin{equation}
\forall m \in \mathbb{N}^*, \ \forall n \in \mathbb{N}, \ \ \ \mathbf{c}_m(n) =\mu\left(\frac{m}{(m,n)}\right) \frac {\varphi(m)}{\varphi\left(\frac{m}{(m,n)}\right)}
\end{equation}

	\section{The proof}
\subparagraph{}
Let $f$ be an arithmetic function such as $ f(n) = \sum_{d|n} g(d)$ where $g$ is also an arithmetic function. 
Thanks to the equality (1) we can write $f$ as a sum,
\begin{equation}
\forall n \in \mathbb{N}^*, \ \ \ f(n) = \sum_{k|n} g(k) = \sum_{k=1}^n \frac{g(k)}{k} \sum_{l=1}^k \exp{\left(\frac{2i\pi nl}{k}\right)} 
\end{equation}
\\
Let's rewrite this sum by changing the indexation,
\begin{equation}
f(n) = \sum_{1 \leq a \leq b \leq n \atop (a,b) = 1} \exp{\left(\frac{2i\pi n a}{b}\right)} C\left(\frac{a}{b}\right)
\end{equation}\\
where,
\begin{equation}
C\left(\frac{a}{b}\right) = \sum_{k \in E\left(\frac{a}{b}\right)} \frac{g(k)}{k}
\end{equation}
with,
\begin{align*}
E\left(\frac{a}{b}\right) &= \Big\{ 1 \leq k \leq n, \ \exists l \leq k,\  \frac{l}{k}=\frac{a}{b} \Big\} \\
&= \Big\{ 1 \leq k \leq n, \ \exists l \leq k, \exists u \in \mathbb{N}^* \vert \ \  (k,l)= (ub,ua) \Big\} \\
&= \Big\{  k = ub \ for \ 1\leq  u \leq \left\lfloor \frac{n}{b} \right\rfloor \Big\}
\end{align*}
Thus when we replace in (6) we have,
\begin{equation}
C\left(\frac{a}{b}\right) = \sum_{u=1}^{\left\lfloor \frac{n}{b} \right\rfloor} \frac{g(ub)}{ub} 
\end{equation}\\
\\
When we use the equality (7) in (5) it comes,\\
\begin{align}
\forall n \in \mathbb{N}^*, \ \ f(n) &= \sum_{1 \leq a \leq b \leq n \atop (a,b) = 1} \exp{\left(\frac{2i\pi n a}{b}\right)}   \sum_{u=1}^{\left\lfloor \frac{n}{b} \right\rfloor} \frac{g(bu)}{bu}  \nonumber\\
&= \sum_{b = 1}^n \sum_{u=1}^{\left\lfloor \frac{n}{b} \right\rfloor} \frac{g(bu)}{bu} \sum_{a=1 \atop (a,b) = 1}^n  \exp{\left(\frac{2i\pi n a}{b}\right)} \nonumber\\
&= \sum_{b = 1}^n \mathbf{c}_b(n)  \sum_{u=1}^{\left\lfloor \frac{n}{b} \right\rfloor} \frac{g(bu)}{bu} 
\end{align}\\
\\
where $\mathbf{c}_b(n)$ is the Ramanujan's sum defined in the preliminaries.\\
\\
Finally, when we replace with the expression (2) given in the preliminaries we have,
\begin{equation}
\forall n \in \mathbb{N}^*, \ \ \ f(n) = \sum_{k=1}^n \mu \left(\frac{k}{(n,k)}\right)  \frac {\varphi(k)}{\varphi\left(\frac{k}{(n,k)}\right)}   \sum_{l=1}^{\left\lfloor\frac{n}{k}\right\rfloor} \frac{g(kl)}{kl}
\end{equation}
which is the identity we wanted to prove.

\section{Application to well known arithmetical functions}

\subsection{Generalized sum of divisors function}
\subparagraph{}
Here we study the function $ \sigma_\gamma(n) = \sum_{d|n} d^\gamma = \sum_{d|n} g(d) $ where $g(d) = d^\gamma$ . This is a generalization of the sum of divisors function ($\gamma =1$) and the number of divisors function ($\gamma = 0$). When we apply our theorem it comes,
\begin {align}
\forall n \in \mathbb{N}^*, \ \ \ \sigma_\gamma(n) &= \sum_{k=1}^n \mathbf{c}_k(n) k^{\gamma-1} \sum_{l=1}^{\left\lfloor\frac{n}{k}\right\rfloor} l^{\gamma-1}\nonumber\\
&=\sum_{k=1}^n \mu \left(\frac{k}{(n,k)}\right)  \frac {\varphi(k)}{\varphi\left(\frac{k}{(n,k)}\right)} k^{\gamma-1} \sum_{l=1}^{\left\lfloor\frac{n}{k}\right\rfloor} l^{\gamma-1}
\end{align}
Let's remark that this expression is true for $\gamma \in \mathbb{R}$ and not only for integers.\\
\subparagraph{•}
For $\gamma = 0$ in (10) we have,

\begin {align}
\forall n \in \mathbb{N}^*, \ \ \ d(n) &= \sum_{k=1}^n \frac{\mathbf{c}_k(n)}{k} H_{\left\lfloor \frac{n}{k}\right\rfloor} \nonumber\\
&= \sum_{k=1}^n \frac{1}{k} H_{\left\lfloor \frac{n}{k} \right\rfloor} \mu\left(\frac{k}{(n,k)}\right) \frac {\varphi(k)}{\varphi\left(\frac{k}{(n,k)}\right)}
\end{align}

where $H_n = \sum_{i=1}^n \frac{1}{i}$ for $n \in \mathbb{N}^*$.\\
\subparagraph{•}
For $\gamma = 1$ in (10) we have,

\begin {align}
\forall n \in \mathbb{N}^*, \ \ \ \sigma(n) &= \sum_{k=1}^n \mathbf{c}_k(n) \left\lfloor \frac{n}{k}\right\rfloor \nonumber\\
&= \sum_{k=1}^n \left\lfloor\frac{n}{k}\right\rfloor  \mu\left(\frac{k}{(n,k)}\right) \frac {\varphi(k)}{\varphi\left(\frac{k}{(n,k)}\right)}
\end{align}

\subsection{A new equality from the product of divisors function}
\subparagraph{}
Let's call $\Pi(n) = \prod_{d|n} d$. Then, we have $log\left(\Pi(n) \right) = \sum_{d|n} log(d)$. We can apply the theorem and it comes,\\
\begin{equation}
\forall n \in \mathbb{N}^*,\ \ \ \log\left(\Pi(n) \right) = \sum_{k=1}^n \mu \left(\frac{k}{(n,k)}\right)  \frac {\varphi(k)}{\varphi\left(\frac{k}{(n,k)}\right)}    \sum_{l=1}^{\left\lfloor\frac{n}{k}\right\rfloor} \frac{log(kl)}{kl}
\end{equation}
Moreover,
\begin{equation*}
\forall n \in \mathbb{N}^*, \ \ \ \Pi^2(n) = \prod_{d|n} d \cdot \prod_{d|n} \frac{n}{d} = \prod_{d|n} n = n^{d(n)}
\end{equation*}
thus, 
\begin{equation}
\forall n \in \mathbb{N}^*, \ \ \ \log\left(\Pi(n) \right) = \frac{1}{2}d(n) \cdot log(n)
\end{equation}
when we equalize (14) and (13) with equality (11) it comes,
\begin{equation}
\forall n \in \mathbb{N}^*, \ \ \ \sum_{k=1}^n \mu \left(\frac{k}{(n,k)}\right)  \frac {\varphi(k)}{\varphi\left(\frac{k}{(n,k)}\right)}  \sum_{l=1}^{\left\lfloor\frac{n}{k}\right\rfloor}  \frac{log \left( \frac{(kl)^2}{n}\right)}{kl}  = 0
\end{equation}
\subsection{Kronecker function}
Lets define Kronecker function as $\delta(n) = \sum_{d|n} \mu(d)$. According to (2) we know that,
\begin{equation*}
\delta(n) =
    \begin{cases}
    1 & if\ \ n = 1 \\
    0 & otherwise
    \end{cases}
\end{equation*} 
Let's apply our theorem to this function,
\begin{equation}
\sum_{k=1}^n \mu \left(\frac{k}{(n,k)}\right)  \frac {\varphi(k)}{\varphi\left(\frac{k}{(n,k)}\right)} \sum_{l=1}^{\left\lfloor\frac{n}{k}\right\rfloor}  \frac{\mu(kl)}{kl} =
    \begin{cases}
    1 & if\ \ n = 1 \\
    0 & otherwise
    \end{cases}
\end{equation}

\end{document}